\documentclass{article}

\def\ra{\rightarrow}

\def\ss{\subseteq}

\def\d{\delta}

\def\P{{\cal P}}
\def\L{{\cal L}}

\def\Re{\hbox{\rm Re}\,}
\def\Im{\hbox{\rm Im}\,}

\def\zbar{\overline{z}}
\def\zetabar{\overline{\zeta}}
\def\O{\Omega}

 \def\HollowBox #1#2{{\dimen0=#1 \advance\dimen0 by -#2       
       \dimen1=#1 \advance\dimen1 by #2                       
        \vrule height #1 depth #2 width #2                    
        \vrule height 0pt depth #2 width #1                   
        \llap{\vrule height #1 depth -\dimen0 width \dimen1}%
       \hskip -#2                                             
       \vrule height #1 depth #2 width #2}}                   
 \def\BoxOpTwo{\mathord{\HollowBox{6pt}{.4pt}}\;}             

\def\endpf{\hfill $\BoxOpTwo$}

\font\teneufm=eufm10
\font\seveneufm=eufm7
\font\fiveeufm=eufm5
\newfam\eufmfam
\textfont\eufmfam=\teneufm
\scriptfont\eufmfam=\seveneufm
\scriptscriptfont\eufmfam=\fiveeufm

\newfam\msbfam
\font\tenmsb=msbm10   \textfont\msbfam=\tenmsb
\font\sevenmsb=msbm7  \scriptfont\msbfam=\sevenmsb
\font\fivemsb=msbm5   \scriptscriptfont\msbfam=\fivemsb
\def\Bbb{\fam\msbfam \tenmsb}

\def\RR{{\Bbb R}}
\def\CC{{\Bbb C}}

\newtheorem{theorem}{Theorem}[section]

\newtheorem{proposition}[theorem]{Proposition}
\newtheorem{lemma}[theorem]{Lemma}

\newtheorem{definition}{Definition}[section]
\newtheorem{example}[definition]{EXAMPLE}

\usepackage{graphicx}

\usepackage{amsmath}

\begin{document}

\begin{center}
{\Large \bf A Tale of Three Kernels}
\bigskip \\
Steven G. Krantz\footnote{Supported in part by a grant from
the National Science Foundation and a grant from the Dean of
Graduate Studies at Washington University.}
\end{center}
\vspace*{.25in}

\begin{quote}
{\bf Abstract:}  \sl
We study three canonical kernels on domains in $\CC$ and $\CC^n$.
We exposit both the history and the substance of these concepts, and
we also present some new calculations and ideas adherent thereto.
\end{quote}

\section{Introduction}

Throughout this paper, a {\it domain} will be a connected open set.
We are interested in three integration kernels that live on
virtually any bounded domain in any complex space.  They are of particular
interest because they are {\it canonical}---they transform
in a natural way under the morphisms in the subject.  They preserve
relevant data and create important new data.  The paradigm for creating
these kernels has broad interest (see [ARO]), and value in many different
subject areas.  We wish to acquaint the readership with this
important set of tools.

This paper is partly expository.  But it also contains a number
of calculations and results that are new, or that at least
are not readily found elsewhere.  We hope that the paper
will be a valuable resource for those who wish to investigate
further in this line of research.

We begin our studies here on domains in the complex plane.  Later in the
paper we shall segue to the function theory of several complex variables.

It is a pleasure to thank Gerald B. Folland and C. Robin Graham for helpful conversations about
the subject matter of this paper.

\section{Basic Definitions}

Let $\Omega$ be a bounded domain in $\CC$.  Define
$$
A^2(\Omega) = \{f \ \hbox{holomorphic on} \ \Omega: \int_\Omega |f(\zeta)|^2 \, dA(\zeta)^{1/2} \equiv
         \|f\|_{A^2(\Omega)} < \infty\} \, .
$$
This is the {\it Bergman space} on $\Omega$.  The space $A^2(\Omega)$ is an inner
product space, in fact a subspace of $L^2(\Omega$), and elementary
arguments (see [KRA2]) show that it is a Hilbert space.	  In particular,
it is complete.

\begin{lemma}   
 Let $K \ss \Omega$ be compact.  There is a constant 
$C_K > 0,$ depending on $K$ and on $n,$ such that

$$ 
\sup_{z \in K} |f(z)| \leq C_K \|f\|_{A^2(\Omega)}\ \ , \ \ \mbox{\rm all}\ f \in A^2(\Omega) . 
$$
\end{lemma}
{\bf Proof:}  Since $K$ is compact, there is an $r = r(K) > 0$ so that, for
any $z \in K$, it holds that $B(z,r) \ss \Omega$.  Therefore, for each $z \in K$ and
$f \in A^2(\Omega)$, we see that (letting $dV$ denote standard Euclidean volume
measure)
\begin{align}
   |f(z)| & =  \frac{1}{V(B(z,r))} \left | \int_{B(z,r)} f(t) \, dV(t) \right | \notag \\
          & \leq  (V(B(z,r)))^{-1/2} \|f\|_{L^2(B(z,r))} \notag \\
          & \leq  c(n) r^{-n} \|f\|_{A^2(\Omega)}  \notag \\
          & \equiv  C_K \| f\|_{A^2(\Omega)} \, .   \tag*{$\BoxOpTwo$}
\end{align}
\vspace*{.15in}

Fix a point $z \in \Omega$.  The functional
$$
e_z:  A^2(\Omega) \ni f \longmapsto f(z)
$$
is easily seen to be a bounded linear functional (take $K$ to be the
singleton $\{z\}$ in the lemma).
Therefore, by the Riesz representation theorem, there is an element
$\varphi_z \in A^2(\Omega)$ such that
$$
f(z) = e_z(f) = \langle f, \varphi_z \rangle
$$
for every $f \in A^2(\Omega)$.  We set $K(z, \zeta) = \overline{\mathstrut \varphi_z(\zeta)}$ and
thus we can write
$$
f(z) = \int_\Omega f(\zeta) K(z,\zeta) \, dA(\zeta) \, ,
$$
where $dA$ is the usual area measure in the plane.  We call $K(z, \zeta)$ the {\it Bergman kernel}
for $\Omega$.  

There is a similar construction\footnote{In
fact Nachman Aronszajn has created an axiomatic theory for reproducing
kernels of this kind.  See [ARO].} for functions on the boundary $\partial \Omega$ of $\Omega \ss \CC^2$.
Assume now, for simplicity, that $\Omega$ has $C^2$ boundary.
Define
$$
H^2(\Omega) = \{f \ \hbox{holomorphic on} \ \Omega: |f|^2 \ \hbox{has a harmonic majorant on} \ \Omega\} \, .
$$
Then it is known (see, for instance, [KRA1]) that this definition of $H^2$ is equivalent
to several other standard and natural definitions of the space.  In particular, if $z \in \Omega$
is fixed, then the functional
$$
\eta_z: H^2(\Omega) \ni f \longmapsto f(z)
$$
is a bounded linear functional.  As a consequence, the Riesz representation theorem tells
us that there is a function $\psi_z \in H^2(\Omega$ such that, for each $f \in H^2(\Omega)$,
$$
f(z) = \eta_z(f) = \langle f, \psi_z \rangle \, .
$$
We set $S(z, \zeta) = \overline{\mathstrut \psi_z(\zeta)}$ and we write
$$
f(z) = \int_{\partial \Omega} f(\zeta) S(z, \zeta) \, d\sigma(\zeta) \, ,
$$
where $d\sigma$ is arc length (or Hausdorff) measure on $\partial \Omega$.  We call $S(z, \zeta)$
the {\it Szeg\H{o} kernel} for $\Omega$.

We note here that most of our results about the Bergman kernel have analogues for the Szeg\H{o} kernel,
with just the same formal proofs.  We shall not provide the details, but encourage the reader
to explore these ideas---or again refer to [ARO] for a broader perspective.

Before we introduce the last of our three kernels, we need some auxiliary information
about the Bergman and Szeg\H{o} kernels.  

\begin{proposition}  \sl 
The Bergman kernel $K(z,\zeta)$ is conjugate symmetric:  $K(z,\zeta) = \overline{K(\zeta,z)}.$
\end{proposition}
{\bf Proof:}  By its very definition, $\overline{K(\zeta,\cdot)} \in A^2(\Omega)$
for each fixed $\zeta.$  Therefore the reproducing property of the
Bergman kernel gives
$$ 
\int_\Omega K(z,t) \overline{K(\zeta,t)} \, dV(t) = \overline{K(\zeta,z)} . 
$$
On the other hand,
\begin{align}
\int_\Omega K(z,t) \overline{K(\zeta,t)} \, dV(t) & =  
                 \overline{\int K(\zeta,t) \overline{K(z,t)} \, dV(t)} \notag \\
                 & =  \overline{\overline{K(z,\zeta)}} = K(z,\zeta) . \tag*{$\BoxOpTwo$}
\end{align}

\begin{proposition}   \sl
 The Bergman kernel is uniquely determined by the
properties that it is an element of $A^2(\Omega)$ in $z,$ is
conjugate symmetric, and reproduces $A^2(\Omega).$
\end{proposition}
{\bf Proof:}  Let $K'(z,\zeta)$ be another such kernel.
Then
\begin{align} 
     K(z,\zeta) & =  \overline{K(\zeta,z)} \notag \\
                & =  \int K'(z,t) \overline{K(\zeta,t)} \, dV(t) \notag \\
                & =  \overline{\int K(\zeta,t) \overline{K'(z,t)} \, dV(t)} \notag \\
                & =  \overline{\overline{K'(z,\zeta)}} \notag \\
		& =  K'(z,\zeta) \, . \tag*{$\BoxOpTwo$}
\end{align} 

Since $L^2(\Omega)$ is a separable Hilbert space then so is its subspace
${A}^2(\Omega).$  Thus there is a complete orthonormal basis $\{\phi_j\}_{j=1}^\infty$
for $A^2(\Omega).$

\begin{proposition}    \sl 
Let $L$ be a compact subset of $\Omega.$  Then the series 
$$
\sum_{j=1}^\infty \phi_j(z) \overline{\phi_j(\zeta)} 
$$ 
sums uniformly
on $L \times L$ to the Bergman kernel $K(z,\zeta).$
\end{proposition}
{\bf Proof:}  By the Riesz-Fischer and Riesz representation theorems, we
obtain
\begin{align}
     \sup_{z \in L} \left ( \sum_{j=1}^\infty |\phi_j(z)|^2 \right )^{1/2} & = 
            \sup_{z \in L} \left \| \{\phi_j(z)\}_{j=1}^\infty\right \|_{\ell^2} \notag \\
             & = \sup_{\scriptsize \|\{a_j\}\|_{\ell^2}=1 
                          \atop
                         \scriptsize z \in L} 
                   \left | \sum_{j=1}^\infty a_j \phi_j(z) \right | \notag \\
             & = \sup_{\scriptsize \|f\|_{A^2}=1 
                          \atop
                             z \in L} 
                  |f(z)| \leq C_L .   \tag{2.4.1}  
\end{align}
In the last inequality we have used Lemma 2.1.
Therefore
$$ 
\sum_{j=1}^\infty \left |\phi_j(z) \overline{\phi_j(\zeta)} \right |
              \leq \left ( \sum_{j=1}^\infty |\phi_j(z)|^2\right )^{1/2}
                   \left ( \sum_{j=1}^\infty |\phi_j(\zeta)|^2\right )^{1/2} 
$$
and the convergence is uniform over $z,\zeta \in L.$  For fixed $z \in \Omega,$
(2.4.1) shows that $\{\phi_j(z)\}_{j=1}^\infty \in \ell^2.$  Hence we have
(for fixed $z$) that $\sum \phi_j(z) \overline{\phi_j(\zeta)} \in \overline{A^2(\Omega)}$ as a function
of $\zeta.$  Let the sum of the series be denoted by ${K}'(z,\zeta).$
Notice that $K'$ is conjugate symmetric by its very definition.  Also,
for $f \in A^2(\Omega),$ we have
$$ 
\int K'(\cdot,\zeta) f(\zeta) \, dV(\zeta) = \sum \widehat{f}(j) \phi_j(\cdot) = f(\cdot), 
$$
where convergence is in the Hilbert space topology.  
[Here $\widehat{f}(j)$ is the $j^{\footnotesize \rm th}$ Fourier coefficient
of $f$ with respect to the basis $\{\phi_j\}.$] 
But Hilbert space convergence
dominates pointwise convergence (Lemma 2.1) so
$$ 
f(z) = \int K'(z,\zeta) f(\zeta) \, dV(\zeta) , \ \ \mbox{\rm all}\ f \in A^2(\Omega) . 
$$
Therefore $K'$ is the Bergman kernel. 
\endpf \smallskip \\ 			  
\medskip \\
\paragraph{Remark:}  It is worth noting explicitly that the proof of Proposition 2.4 shows that
$$
\sum \phi_j(z) \overline{\phi_j(\zeta)}
$$ 
equals the Bergman kernel $K(z,\zeta)$
{\em no matter what the choice} of complete orthonormal basis $\{\phi_j\}$
for $A^2(\Omega).$
\smallskip   \\

\begin{proposition}    \sl
If $\Omega$ is a bounded domain in $\CC$ then the mapping
$$ 
 P: f \mapsto \int_{\Omega} K(\cdot,\zeta) f(\zeta) \, dV(\zeta) 
$$
is the Hilbert space orthogonal 
projection of $L^2(\Omega,\, dV)$ onto $A^2(\Omega).$
\end{proposition}
{\bf Proof:}  Notice that $P$ is idempotent and self-adjoint and that
$A^2(\Omega)$ is precisely the set of elements of $L^2$ that are fixed
by $P.$ 
\endpf \smallskip \\

\begin{definition}   \rm
Let $\Omega \ss \CC$ be a domain and let $f: \Omega \ra \CC$ be a 
holomorphic function.  It is sometimes convenient to write $z = z + iy$ 
and $f = u + iv$.  Then we may render $f$ as a {\it real mapping}
$$
(x,y) \longmapsto \bigl ( u(x + iy), v(x + iy) \bigr ) \, .
$$
Thus it makes sense to consider the {\it real Jacobian} matrix
$$
J_\RR f(z) = \left ( \begin{array}{cc}
                    \displaystyle \frac{\partial u}{\partial x} & \displaystyle \frac{\partial u}{\partial y} \\  [.15in]
		    \displaystyle \frac{\partial v}{\partial x} & \displaystyle \frac{\partial v}{\partial y} \\
		     \end{array}
	     \right ) \, .
$$
\end{definition}

It is natural to wonder what is the relationship between the real Jacobian
$J_\RR f$ and the complex derivative $f$.  The next proposition enunciates
the definitive result.

\begin{proposition} \sl   
If $f$ is a holomorphic function on a domain $\Omega \ss \CC$ then
$$
\hbox{det} \ J_\RR f(x,y) = |f'(z)|^2 \, .
$$
\end{proposition}
{\bf Proof:}  This follows immediately from the Cauchy-Riemann equations,
or see [GRK]. 
\endpf 
\smallskip \\

A holomorphic mapping $f: \Omega_1 \ra \Omega_2$ of domains
$\Omega_1 \ss \CC, \Omega_2 \ss \CC$ is said to be {\em
biholomorphic} if it is one-to-one, onto, and $\mbox{\rm
det}\, {J}_\CC f(z) \not = 0$ for every $z \in \Omega_1$. It
is automatic that the inverse mapping is holomorphic (see
[GRK]). In the context of one complex variable, such a mapping
is often called {\it conformal}.

In what follows we denote the Bergman kernel for a given domain
$\Omega$ by $K_\Omega.$

\begin{proposition}   \sl 
Let $\Omega_1, \Omega_2$ be domains in $\CC.$
Let $f: \Omega_1 \ra \Omega_2$ be biholomorphic.  Then
$$
f'(z)  K_{\Omega_2}(f(z),f(\zeta)) \overline{f'(\zeta)}
                         = K_{\Omega_1}(z,\zeta) \, . 
$$
\end{proposition}
{\bf Proof:} Let $\phi \in A^2(\Omega_1).$ Then, by change of variable,
\begin{eqnarray*}
\lefteqn{\int_{\Omega_1}  
        f'(z)  K_{\Omega_2}(f(z),f(\zeta))  \overline{f'(\zeta)} \phi(\zeta) \, dV(\zeta)} \\
& = & \int_{\Omega_2}   f'(z)  K_{\Omega_2}(f(z),\widetilde{\zeta}) 
            \overline{f'(f^{-1}(\widetilde{\zeta}))}
               \phi(f^{-1}(\widetilde{\zeta})) \\
& \mbox{} & \qquad \times  \hbox{det}\, J_\RR \bigl [f^{-1} \bigr ] (\widetilde{\zeta}) \, dV(\widetilde{\zeta}) .
\end{eqnarray*}
By Proposition 2.6 this simplifies to 
$$  
  f'(z) \int_{\Omega_2} K_{\Omega_2}(f(z),\widetilde{\zeta}) 
         \left \{ \left (  f'(f^{-1}(\widetilde{\zeta})) \right )^{-1} \phi\left ( f^{-1}(\widetilde{\zeta}) \right )\right \}  \, dV(\widetilde{\zeta}) . 
$$
By change of variables, the expression in braces $\{ \mbox{\ \ } \}$ is an 
element of ${A}^2(\Omega_2).$  So the reproducing property of 
$K_{\Omega_2}$ applies and the last line equals
$$  
=    f'(z) \left [ f'(z) \right ]^{-1}
               \phi\left (f^{-1}(f(z))\right ) = \phi(z) . 
$$
By the uniqueness of the Bergman kernel, the proposition follows.
\endpf \smallskip \\

\begin{proposition}   \sl
 For $z \in \Omega \subset \subset \CC$
it holds that $K_\Omega(z,z) > 0.$
\end{proposition}
{\bf Proof:}  Now
$$ 
K_\Omega(z,z) = \sum_{j=1}^\infty |\phi_j(z)|^2 \geq 0 .
$$
If in fact $K(z,z) = 0$ for some $z$ then $\phi_j(z) = 0$ for all
$j$ hence $f(z) = 0$ for every $f \in A^2(\Omega).$  This is absurd.
\endpf \smallskip \\

\begin{proposition}  \sl  
Let $\Omega \subset \subset \CC$ be a domain.  Let $z \in \Omega.$
Then
$$  
K(z,z) = \sup_{f \in A^2(\Omega)} \frac{|f(z)|^2}{\|f\|^2_{A^2}}
                = \sup_{\|f\|_{A^2(\Omega)}=1} |f(z)|^2 . 
$$
\end{proposition}
{\bf Proof:}  Now
\begin{eqnarray*}
   K(z,z) & = & \sum | \phi_j(z) |^2 \\
          & = & \left ( \sup_{\|\{a_j\}\|_{\ell^2}=1} \left | \sum \phi_j(z) a_j \right | \right )^2 \\
          & = & \sup_{\|f\|_{A^2} = 1} |f(z)|^2 , 
\end{eqnarray*}
by the Riesz-Fischer theorem,
$$   
= \sup_{f \in A^2} \frac{|f(z)|^2}{\|f\|_{A^2}^2} .  \eqno \BoxOpTwo 
$$
\vspace*{.12in}

There is in fact a third kernel that will be of interest for us
here. It was discovered fairly recently by Lu Qi-Keng Hua (see
[HUA] and also [KOR]). And it is not as well known as it
should be.

If $\Omega \ss \CC$ us a bounded domain with $C^2$ boundary then let $S(z, \zeta)$ be
its Szeg\H{o} kernel.  We define the {\it Poisson-Szeg\H{o}} kernel of $\Omega$ to
be 
$$
{\cal P}(z, \zeta) = \frac{|S(z,\zeta)|^2}{S(z,z)} \, .
$$
Note that $S(z,z) \ne 0$ (by a proof analogous to that for Proposition 2.8) 
so this last definition makes sense.  Now we have the 
following fundamental result:

\begin{proposition} \sl
Let $f$ be a continuous function on $\overline{\Omega}$ which is holomorphic
on $\Omega$.  Then, for any $z \in \Omega$,
$$
f(z) = \int_{\partial \Omega} f(\zeta) {\cal P}(z, \zeta) \, d\sigma(\zeta) \, .
$$
\end{proposition}					     
{\bf Proof:}  Fix $z \in \Omega$.  Define $g(\zeta) = \overline{S(z,\zeta)} \cdot f(\zeta)/S(z,z)$.
Then it is easy to see that $g \in H^2(\Omega)$ as a function of $\zeta$.  As a result,
\begin{align}
\int_{\partial \Omega} f(\zeta) {\cal P}(z, \zeta) \, d\sigma(\zeta) & = 
         \int_{\partial \Omega} \left [ f(\zeta) \cdot \frac{\overline{S(z,\zeta)}}{S(z,z)} \right ] \cdot S(z, \zeta) \, d\sigma(\zeta) \notag \\
             & =  \int_{\partial \Omega} g(\zeta) \cdot S(z, \zeta) \, d\sigma(\zeta) \notag \\
	     & =  g(z)  \notag \\
	     & =  f(z) \, .  \tag*{$\BoxOpTwo$}	  \\ \notag
\end{align}

The Poisson-Szeg\H{o} kernel is, by the last proposition, a reproducing kernel.  But it is also positive,
which is a very useful property (as we often see with the classical Poisson kernel).  	See also
Section 8 where this idea is used decisively.

\section{Calculating the Kernels}

It is in general rather difficult to explicitly calculate any of the three kernels
being discussed here.  [A similar statement can be made about the classical
Poisson kernel, but see for instance [KRA3].]	This situation is in marked contrast
to that for the Cauchy kernel.  That kernel is the same for {\it every domain}, and
is explicitly given by
$$
\frac{1}{2\pi i} \, \frac{1}{\zeta - z} \, .
$$
The Bergman, Szeg\H{o}, Poisson-Szeg\H{o}, and Poisson kernels are typically different for different 
domains.  Unless the domain has a good deal of symmetry (like the unit ball, for example), there
is little hope of actually writing down a formula for one of these kernels.

For the moment we shall content ourselves with writing down the Bergman and Szeg\H{o} kernels
for the disc.  Although there are several ways to do this (see [KRA2]), we shall make
good use now of Proposition 2.4.  

\begin{example} \rm
Our domain now is $\Omega = D$ the unit disc.  Consider the functions
$$
\psi_j(z) = z^j \ , \qquad j = 0, 1, 2, \dots 
$$
on $\Omega$.  By parity, these functions are orthogonal in the $L^2(\Omega)$ inner product.  By
standard results on power series of holomorphic functions, the set $\{\psi_j\}$ forms
a complete orthogonal system in $A^2(\Omega)$.  We calculate that
\begin{eqnarray*}
\|\psi_j\|_{A^2(\Omega)}^2 & = & \int_0^{2\pi} \int_0^1 |r e^{i \theta}|^{2j} r \, dr d\theta   \\
			   & = & 2\pi \int_0^1 r^{2j+1} \, dr \\
			   & = & \frac{\pi}{j+1} \, .
\end{eqnarray*}
Thus
$$
\|\psi_j\|_{A^2(\Omega)} = \sqrt{\frac{\pi}{j+1}}
$$
and
$$
\varphi_j(z) \equiv \sqrt{\frac{j+1}{\pi}} \, z^j
$$
is a complete orthonormal system on $\Omega$.  Thus Proposition 2.4 tells us that
the Bergman kernel on $\Omega$ is given by
$$
K(z,\zeta) = \sum_{j=1}^\infty \frac{j+1}{\pi} z^j \overline{\zeta}^j \, .  \eqno (3.1.1)
$$

We know that 
$$
\sum_{j=0}^\infty \lambda^j = \frac{1}{1 - \lambda}
$$
when $|\lambda| < 1$.  Rewriting this as
$$
\sum_{j=-1}^\infty \lambda^{j+1} = \frac{1}{1 - \lambda}
$$
and differentiating in $\lambda$, we obtain
$$
\sum_{j=0}^\infty (j+1) \lambda^j = \frac{1}{(1 - \lambda)^2} \, .
$$
Applying this last formula to $(3.1.1)$ yields
$$
K(z,\zeta) = \frac{1}{\pi} \frac{1}{(1 - z\cdot \overline{\zeta})^2} \, .
$$
That is the Bergman kernel for the disc.
\end{example}

\begin{example}  \rm
Our domain once again is $\Omega = D$ the unit disc.  Consider the functions
$$
\psi_j(z) = z^j \ , \qquad j = 0, 1, 2, \dots 
$$
on $\Omega$.  By parity, these functions are orthogonal in the $L^2(\partial \Omega)$ inner product.  By
standard results on power series of holomorphic functions, the set $\{\psi_j\}$ forms
a complete orthogonal system in $H^2(\Omega)$.  We calculate that
\begin{eqnarray*}
\|\psi_j\|_{H^2(\Omega)}^2 & = & \int_0^{2\pi} |e^{i \theta}|^{2j} d\theta   \\
			   & = & 2\pi 
\end{eqnarray*}
Hence
$$
\varphi_j(z) \equiv \frac{1}{\sqrt{2\pi}} \, z^j
$$
is a complete orthonormal system in $H^2(\Omega)$.  Applying a suitable variant of Proposition 2.4 (adapted
to $H^2$ and the Szeg\H{o} kernel), we find that
$$
S(z,\zeta) = \sum_{j=0}^\infty \frac{1}{2\pi} z^j \overline{\zeta}^j = \frac{1}{2\pi} \, \frac{1}{1 - z \cdot \overline{\zeta}} \, .
$$
This is the Szeg\H{o} kernel for the unit disc.
\end{example}

Observe that the Szeg\H{o} kernel on the disc has a form similar to the Bergman kernel on the disc, but the singularity
is one degree lower.  This makes sense, because the Szeg\H{o} integral is a boundary integral while the Bergman
integral is an integral over the 2-dimensional domain.  We shall say more about these relationships in the
next section.

\begin{example} \rm
Our domain again is $\Omega = D$ the unit disc.   Let us calculate the Poisson-Szeg\H{o}
kernel for this domain.

Now
$$
{\cal P}(z, \zeta) = \frac{|S(z,\zeta)|^2}{S(z,z)} = \frac{1}{2\pi} \, \frac{1/|1 - z \cdot \overline{\zeta}|^2}{1/[1 - |z|^2]}
		   = \frac{1}{2\pi} \, \frac{1 - |z|^2}{|1 - z \cdot \overline{\zeta}|^2} \, .
$$
This is the Poisson-Szeg\H{o} kernel for the disc.  At this point it is interesting to introduce polar coordinates:
We set $z = r e^{i\theta}$ and $\zeta = e^{i\psi}$.  Then we find that
$$
{\cal P}(r e^{i\theta}, e^{i\psi} ) = \frac{1}{2\pi} \frac{1 - r^2}{\left |1 - r e^{i(\theta - \psi)}\right |^2} \, .
$$
We have discovered that the Poisson-Szeg\H{o} kernel on the disc is actually the classical Poisson kernel!
\end{example}

\begin{example} \rm
Let us endeavor to calculate the Bergman kernel for the planar annulus
$$
A = \{z \in \CC: 1 < |z| < 2\} \, .
$$
As Bergman himself points out in [BER], such a calculation is
essentially intractable---it would involve elliptic functions.
But we can make some interesting qualitative statements about
the kernel.

It is easy to see that the functions
$$
\psi_j(z) = z^j \ , \qquad j = \dots, -3, -2, -1, 0, 1, 2, 3, \dots
$$
form a complete orthogonal system on $A$.  Moreover,
\begin{eqnarray*}
\|\psi\|_{A^2(A)}^2 & = & \int_0^{2\pi} \int_1^2 |r e^{i\theta}|^{2j} \cdot r \, dr d\theta  \\
		    & = &  2\pi \int_1^2 r^{2j+1} \, dr \\   [.2in]
		    & = & \left \{ \begin{array}{lcr}
		              \frac{2\pi}{2j+2} \left [ 2^{2j+2} - 1 \right ] & \hbox{if} & j \ne -1 \\  [.15in]
			      2\pi \, \log 2                                       & \hbox{if} & j = - 1 \, .
				   \end{array}
			  \right.
\end{eqnarray*}
Thus we find that
$$
\|\psi\|_{A^2(A)} = \left \{ \begin{array}{lcr}
                        \sqrt{\frac{2\pi}{2j+2}} \, \sqrt{2^{2j+2} - 1} & \hbox{if} & j \ne - 1 \\	[.15in]
			     \sqrt{2\pi \, \ln 2} & \hbox{if} & j = - 1 \, .
			     \end{array}
		    \right.
$$

As usual, we conclude (using Proposition 2.4) that
\begin{eqnarray*}
K_A(z, \zeta) & = & \sum_{j=-\infty \atop
                        j \ne -1}^\infty \frac{j+1}{\pi(2^{2j+2} -1)} z^j \overline{\zeta}^j + \frac{1}{2\pi \ln 2} z^{-1} \overline{\zeta}^{-1} \\
              & = & \sum_{j\leq -2} + \sum_{j=-1} + \sum_{j\geq0} \\  [.12in]
	      & \equiv & I + II + III \, .
\end{eqnarray*}
We shall analyze $I$, $II$, and $III$ separately.

In fact $II$ is of little interest.  It is bounded on the annulus, and all its derivatives are bounded. It
induces a bounded operator on any $L^p$ or Sobolev space, and in fact it is a smoothing operator
in the sense of pseudodifferential operators (see [KRA4]).  So it is trivial from our point of view.

Term $I$ is more interesting.  We write
\begin{eqnarray*} 
I & = & \sum_{j=-\infty}^{-2} \frac{j+1}{\pi(2^{2j+2} - 1)} z^j \zetabar^j \\
  & = & \sum_{j=-\infty}^{-2} - \frac{j+1}{\pi} z^j \zetabar^j + \sum_{j=-\infty}^{-2} \left [
                        \frac{j+1}{\pi(2^{2j+2} - 1)} + \frac{j+1}{\pi} \right ] z^j \zetabar^j \\
  & = & \sum_{j=-\infty}^{-2} - \frac{j+1}{\pi} z^j \zetabar^j + \sum_{j=-\infty}^{-2} \frac{j+1}{\pi}
                        \frac{2^{2j+2}}{2^{2j+2} - 1} z^j \zetabar^j \\
  & \equiv & I_1 + I_2 \, .
\end{eqnarray*}
Recalling that $|z| \leq 1$, $|\zeta| \leq 1$ on the annulus (and also keeping in mind
that $j \leq -2 < 0$), we see that $I_2$ is bounded, and all its
derivatives are bounded.  So, as we noted above about $II$, this term is trivial from
the point of view of pseudodifferential operators.  We may ignore it.

As for $I_1$, we note that
\begin{eqnarray*}
\sum_{j=-\infty}^{-2} - \frac{j+1}{\pi} \lambda^j & = & \frac{d}{d\lambda} \sum_{j=-\infty}^{-2} - \frac{1}{\pi} \lambda^{j+1} \\
		& = & \frac{d}{d\lambda} \sum_{j=0}^\infty - \frac{1}{\pi} \lambda^{-j-1} \\
		& = & \frac{d}{d\lambda} \left [ - \frac{1}{\pi} \frac{1}{\lambda - 1} \right ] \\
		& = & \frac{1}{\pi} \cdot \frac{1}{(\lambda - 1)^2} \, .
\end{eqnarray*}
This last is valid as long as $|\lambda| > 1$.  We conclude that
$$
I_1 = \frac{1}{\pi} \frac{1}{(1 - z \cdot \zetabar)^2} \, .
$$

It remains to analyze $III$.

Now we write
\begin{eqnarray*}
III & = & \sum_{j=0}^\infty \frac{j+1}{\pi \cdot 2^{2j+2}} z^j \overline{\zeta}^j   \\
   && \quad + \sum_{j=0}^\infty \left [ \frac{j+1}{\pi(2^{2j+2} - 1)} - \frac{j+1}{\pi \cdot 2^{2j+2}} \right ] z^j \overline{\zeta}^j   \\
 & \equiv & III_1(z, \zeta) + III_2(z, \zeta)  \, .
\end{eqnarray*}
					 
But of course
$$
III_2(z, \zeta) = \sum_{j=0}^\infty \frac{\pi(j+1)}{[2^{2j+2} - 1] \cdot 2^{2j+2}} z^j \zetabar^j \, .
$$
Since $|z| < 2$ and $|\zeta| < 2$ on our domain, we see that the series for $E$ converges absolutely
and uniformly; and the same can be said for any derivative of $E$.  As a result, the error term $E$
is trivial from our point of view.

Thus we are left with studying 
$$
III_1(z, \zeta) = \sum_{j=0}^\infty \frac{j+1}{\pi \cdot 2^{2j+2}} z^j \overline{\zeta}^j \, .
$$
As in our study of the Bergman kernel on the disc, we note that
\begin{eqnarray*}
\sum_{j=0}^\infty \frac{j+1}{2^{2j+2}} \cdot \lambda^j & = &
             \frac{1}{4} \, \sum_{j=0}^\infty (j+1) \cdot \left ( \frac{\lambda}{4} \right )^j \\
	     & = & \frac{d}{d\lambda} \sum_{j=0}^\infty \left ( \frac{\lambda}{4} \right )^{j+1} \\
	     & = & \frac{d}{d\lambda} \left ( \frac{\lambda}{4} \cdot \frac{1}{1 - \lambda/4} \right ) \\
	     & = & \frac{4}{(4 - \lambda)^2} \, .
\end{eqnarray*}
We conclude that
$$
III_1(z, \zeta) = \frac{4}{\pi} \cdot \frac{1}{(4 - z\cdot \zetabar)^2} \, .
$$

Putting all of our calculations together, we find that
$$
K_A(z, \zeta) \approx \frac{4}{\pi} \cdot \frac{1}{(4 - z \cdot \zetabar)^2} + \frac{1}{\pi} \cdot \frac{1}{(1 - z \cdot \zetabar)^2} \, .
$$
Here we have omitted the error terms, all of which we have determined to be trivial in this context.

We see that the Bergman kernel for the annulus is, up to negligible error terms, the sum of the Bergman kernel for the
disc of radius 2 and the Bergman kernel for the disc of radius 1 (or the complement of this disc---via Proposition 2.7
and the mapping $z \mapsto 1/z$).  This is quite satisfying, as it is 
consistent with the geometry of the situation.
\end{example}

\section{The Relevance of Stokes's Theorem}

It is natural to wonder how the Cauchy integral formula fits into the ideas we have been discussing here.
We have already seen some connection, because the Poisson-Szeg\H{o} kernel on the disc turns out
to be the classical Poisson kernel (and of course the Poisson kernel is essentially the real
part of the Cauchy kernel---see [KRA5]).   But we can make the connection even more explicit
by utilizing Stokes's theorem.

Because we are doing complex analysis here, it is useful to have a version of Stokes's theorem
that is written in that language.  Recall that if $z = z + iy$ is a complex variable then
$$
\frac{\partial}{\partial z} = \frac{1}{2} \left [ \frac{\partial}{\partial x} - i \frac{\partial}{\partial y} \right ] \qquad
  \hbox{and} \qquad
\frac{\partial}{\partial \zbar} = \frac{1}{2} \left [ \frac{\partial}{\partial x} + i \frac{\partial}{\partial y} \right ] \, .
$$
It follows that
$$
\begin{array}{cc}
\displaystyle \frac{\partial}{\partial z} \, z = 1 \quad & \quad \displaystyle \frac{\partial}{\partial z} \, \zbar = 0 \, , \\ [.15in]
\displaystyle \frac{\partial}{\partial \zbar} \, z = 0 \quad & \quad \displaystyle \frac{\partial}{\partial \zbar} \, \zbar = 1 \, . \\
\end{array}
$$

We also define
$$
dz = dx + i dy \qquad \hbox{and} \qquad d\zbar = dx - i dy \, .
$$
The differentials $dz$ and $\zbar$ pair with the tangent vectors $\partial/\partial z$ and
$\partial/\partial \zbar$ in the expected manner.

Now we have

\begin{theorem}[Stokes]  \sl
Let $\Omega \ss \CC$ be a bounded domain with $C^1$ boundary.  Let $f$ be a function
that is $C^1$ on $\overline{\Omega}$.  Then
$$
\oint_{\partial \Omega} f(z) \, dz = \mathop{\int \!\!\! \int}_\Omega \, \frac{\partial f}{\partial \zbar} \ d\zbar \wedge dz \, . \eqno (4.1.1)
$$
\end{theorem}

\noindent We invite the reader to write out the formula (4.1.1) in real coordinates so as to
be able to recognize it as the Stokes's theorem that can be found in any calculus book [BLK].

Now let us apply Theorem 4.1 to learn something about the Bergman kernel.  Let $\Omega = D$, the
unit disc.  Let $f$ be a function that is $C^1$ on $\overline{D}$ and holomorphic on $D$.  
Citing the Cauchy integral formula, we write (for $z \in D$, $\zeta \in D$, $\zeta = \xi + i\eta$)
\begin{eqnarray*}
  f(z) & = & \frac{1}{2\pi i} \oint_{\partial D} \, \frac{f(\zeta)}{\zeta - z} \, d\zeta \\
       & = & \frac{1}{2\pi i} \oint_{\partial D} \, \frac{f(\zeta) \cdot \zetabar}{1 - z \cdot \zetabar} \, d\zeta \\
       & = & \frac{1}{2\pi i} \mathop{\int\!\!\!\int}_D \, \frac{\partial}{\partial \zetabar} \left [
			 \frac{\zetabar}{1 - z \cdot \zetabar} \right ] f(\zeta) \, d\zetabar \wedge d\zeta \\
       & = & \frac{1}{2\pi i} \mathop{\int\!\!\!\int}_D \frac{1}{(1 - z \cdot \zetabar)^2} \cdot f(\zeta) \, 2i \, d\xi d\eta \\
       & = & \frac{1}{\pi} \mathop{\int\!\!\!\int}_D \, \frac{1}{(1 - z \cdot \zetabar)^2} \cdot f(\zeta) \, d\xi d\eta  \\
       & = & \mathop{\int\!\!\!\int}_D K_D(z, \zeta) \, f(\zeta) \, d\xi d\eta \, .
\end{eqnarray*}

We see that, beginning with the familiar Cauchy integral formula, we have discovered the Bergman reproducing 
formula.

Now let us think about the Szeg\H{o} reproducing formula.  We write, using the notation
$z = r e^{i\theta}$, $\zeta = e^{i\psi}$, and with $d\sigma$ denoting arc length measure on $\partial D$,
\begin{eqnarray*}
  f(z) & = &  \frac{1}{2\pi i} \oint_{\partial D} \, \frac{f(\zeta)}{\zeta - z} \, d\zeta \\
       & = &  \frac{1}{2\pi i} \int_0^{2\pi} \frac{f(e^{i\psi}}{e^{i\psi} - r e^{i\theta}} \, i e^{i\psi} \, d\psi \\
       & = &  \frac{1}{2\pi}   \int_0^{2\pi} \frac{f(e^{i\psi}}{1 - r e^{i(\theta - \psi)}} \, d\psi \\
       & = &  \frac{1}{2\pi}   \int_{\partial D} \, \frac{f(\zeta)}{1 - z \cdot \zetabar} \ d\sigma(\zeta) \, .
\end{eqnarray*}

We see that, beginning with the classical Cauchy integral formula, we have rediscovered the Szeg\H{o} reproducing
formula.  Since everything takes place on the boundary in this argument, we did not need to use Stokes's theorem this time.

The calculations that we have done in this section suggest that the Bergman kernel is, in general, a derivative
of the Szeg\H{o} kernel.  Certainly for the domain the disc this is the case.  In fact this paradigm is
true in some generality---see [NRSW].  

\section{Boundary Behavior of the Bergman Kernel}

It is a fact that, on a smoothly bounded domain $\Omega \ss \CC$, the Bergman kernel is smooth on
$\overline{\Omega} \times \overline{\Omega} \setminus \bigtriangleup$, where $\bigtriangleup$
is the boundary diagonal $\partial D \times \partial D$.  A fancy way to see this is to note that the Bergman kernel $K_\Omega(z, \zetabar)$
is equal to the Bergman projection $P$ of the Dirac delta mass $\delta_\zeta$ and then 
to cite the pseudolocality of $P$ which can be proved using regularity results of the $\overline{\partial}$-Neumann
problem.  A more prosaic way to prove the assertion, at least in the case of simply connected $\Omega$, 
is to let $\varphi: \Omega \ra D$ be the Riemann mapping function.  By a classical result of
Painlev\'{e} (see [BEK]), the mapping $\varphi$ and its inverse $\varphi^{-1}$ extend to smooth
diffeomorphisms of the respective closures.  Thus the transformation formula
$$
\varphi'(z) \cdot K_D(\varphi(z), \varphi(\zeta)) \cdot \overline{\varphi'(\zeta)} = K_\Omega(z, \zeta)
$$
shows that $K_\Omega$ blows up at the boundary in just the same way as $K_D$ blows up at the 
boundary.  Of course it is clear by inspection that
$$
K_D(z, \zeta) = \frac{1}{\pi} \, \frac{1}{(1 - z \cdot \zetabar)^2}
$$
is smooth on $\overline{D} \times \overline{D} \setminus \bigtriangleup$.  So the same assertion
is true for $\Omega$.

Put in other words, the Bergman kernel $K_\Omega(z, \zeta)$ for a smoothly bounded domain can only blow
up as $z, \zeta \in \Omega$ approach the {\it same} boundary point $P$.  It is worthwhile to
study the rate of blowup of the kernel in such a circumstance.  In fact that rate of blowup can be used
to classify domains (see [APF]).

\begin{example} \rm
Let $\Omega = D$, the unit disc.  Fix $P = 1 + i0 \in \partial D$.  For $\delta > 0$ small,
let $z_\delta = (1 - \delta) + i0$ and $\zeta_\delta = (1 - \delta) + i0$.  Then
$$
K_D(z_\delta, \zeta_\delta) = \frac{1}{\pi} \, \frac{1}{\left (1 - (1 - \delta)\cdot (1 - \delta)\right )^2}
   = \frac{1}{\pi} \, \frac{1}{(2\delta - \delta^2)^2} \approx \frac{1}{4\pi} \, \frac{1}{\delta^2} \, .
$$
Note that the distance of either $z_\delta$ or $\zeta_\delta$ to the boundary of $D$ is precisely
$\delta$.  So we may write
$$
K_D(z_\delta, \zeta_\delta) \approx c \cdot \frac{1}{[\hbox{dist}(z_\delta, \partial D)]^2} 
                              = c \cdot \frac{1}{[\hbox{dist}(\zeta_\delta, \partial D)]^2} \, .
$$
\end{example}

\begin{example} \rm
Now let ${\cal Q} = \{z \in \CC: \Re z > 0, \Im z > 0\}$.  This is the upper-right quarter
plane.  It is of course conformally equivalent to the half plane
$U = \{z \in \CC: \Im z > 0\}$, and that is in turn conformally equivalent to the disc.
We shall use this information to calculate the Bergman kernel of ${\cal Q}$.

Let
\begin{eqnarray*}
\varphi: U & \longrightarrow & D \\
         z & \longmapsto & \frac{i - z}{i + z} \, .
\end{eqnarray*}
It is a straightforward matter to check that $\varphi$ maps $U$ conformally
onto $D$.  Also $\varphi'(z) = -2i/(i + z)^2$.

Now let us apply Proposition 2.7 to determine the Bergman kernel of the
upper halfplane $U$.  We have that
\begin{eqnarray*}
K_U(z, \zeta) & = & \frac{-2i}{(i + z)^2} \cdot \frac{1}{\pi} \cdot
            \frac{1}{\left ( 1 - \left ( \displaystyle \frac{i - z}{i + z} \right ) \cdot
	                         \left ( \displaystyle \frac{- i - \zetabar}{- i + \zetabar} \right ) \right )^2} \cdot
				 \frac{2i}{(-i + \zetabar)^2} \\
	      & = & \frac{4}{\pi} \cdot \frac{1}{\left ( (i + z)(-i + \zetabar) - (i - z)(-i - \zetabar) \right )^2} \\
	      & = & - \frac{1}{\pi} \cdot \frac{1}{( \zetabar - z)^2} \, .
\end{eqnarray*}

\noindent Note that, if we let $z_\delta = x + i \delta$ and $\zeta_\delta = x + i \delta$, then
$$
K_U(z_\delta, \zeta_\delta) = - \frac{1}{\pi} \cdot \frac{1}{\left ( (x - i \delta) - (x + i \delta) \right )^2}
                          = \frac{1}{4\pi} \cdot \frac{1}{\delta^2} \, .
$$
Thus, not surprisingly,
$$
K_U(z_\delta, \zeta_\delta) = c \cdot \frac{1}{\hbox{dist}(z_\delta, \partial U)^2} =
                              c \cdot \frac{1}{\hbox{dist}(\zeta_\delta, \partial U)^2} \, ,
$$
just as in the last example.

Now consider the upper-right quarter plane ${\cal Q}$.  We shall calculate the Bergman kernel for ${\cal Q}$ by again using
Proposition 2.7.  Notice that $\varphi(z) = z^2$ is a conformal mapping of ${\cal Q}$ onto $U$.
Thus we have that
\begin{eqnarray*}
K_{\cal Q}(z, \zeta) & = & 2z \cdot \frac{-1}{\pi} \cdot \frac{1}{(\zetabar^2 - z^2)^2} \cdot 2\zetabar =
	      - \frac{1}{\pi} \cdot \frac{4z \zetabar}{(\zetabar^2 - z^2)^2} \, .
\end{eqnarray*}			      

Now let $z_\delta = \delta + i \delta$ and $\zeta_\delta = \delta + i \delta$.  Then 
$$
K_{\cal Q}(z_\delta, \zeta_\delta) = - \frac{1}{\pi} \cdot \frac{4 (\d + i \d)(\d - i \d)}{\left (
                          (\d - i \d)^2 - (\d + i \d)^2 \right )}
			 = \frac{2}{\pi i} \cdot \frac{1}{\delta^2} \, .
$$
What is interesting now is that the distance of $z_\d$ or $\zeta_\d$ to $\partial {\cal Q}$ is
$\delta$, so that
$$
K_{\cal Q}(z_\d, \zeta_\d) = c \cdot \frac{1}{\hbox{dist}(z_\d, \partial {\cal Q})^2} =
                             c \cdot \frac{1}{\hbox{dist}(\zeta_\d, \partial {\cal Q})^2} \, .
$$
Yet the distance of $z_\delta$ or $\zeta_\d$ to the limit point 0 is $\sqrt{\delta}$.  Thus
$$
K_{\cal Q}(z_\d, \zeta_\d) = c \cdot \frac{1}{\hbox{dist}(z_\d, 0)^4} = c \cdot \frac{1}{\hbox{dist}(\zeta_\d, 0)^4} \, .
$$
\end{example}
This last example gives a way to distinguish a smooth boundary point (as for the disc
$D$ or the upper halfplane $U$) from a boundary point with a corner (as for ${\cal Q}$).

\section{Calculation of the Szeg\H{o} Projection}

The operator
$$
Sf(z) = \int_{\partial \Omega} f(\zeta) S(z,\zeta) \, d\sigma(\zeta)
$$
is the orthogonal Hilbert space projection of $L^2(\partial \O)$ onto $H^2(\O)$.  It
is an important operator for harmonic analysis and complex function theory.  

It is both satisfying and enlightening to calculate some Szeg\H{o} projections of
particular functions.  This
we now do.

\begin{example} \rm
Let $\O = D$, the disc.  Let $f(z) = \overline{z}$.  If $\alpha(z)$ is any
element of $H^2(D)$ then we see that
$$
\langle \alpha, f \rangle = \int_{\partial D} \alpha(z) \cdot \overline{\overline{z}} \, d\sigma(z) 
        = \int_{\partial D} \alpha(z) \cdot z \, d\sigma(z) = 0
$$
by the mean value property of holomorphic functions.  Thus $f$ is orthogonal to the space $H^2$.
We now confirm that assertion with a calculation:
\begin{eqnarray*}
Sf(z) & = & \frac{1}{2\pi} \int_{\partial D} \frac{\zetabar}{(1 - z \cdot \zetabar} \, d\sigma(\zeta) \\
      & = & \frac{1}{2\pi i} \oint_{\partial D} \frac{\zetabar^2}{1 - z \cdot \zetabar} \, d\zeta \, .
\end{eqnarray*}
We use here the fact that $d\zeta = i e^{it} \, dt$.  Now this last
\begin{eqnarray*}
  & = & \frac{1}{2\pi i} \oint_{\partial D} \frac{1}{\zeta(\zeta - z)} \, d\zeta \\
  & = & \frac{1}{2\pi i} \oint_{\partial D} \frac{-1}{z} \left ( \frac{1}{\zeta} - \frac{1}{\zeta - z} \right ) \, d\zeta \\
  & = & \frac{-1}{z} \cdot (1 - 1) \\
  & = & 0 \, .
\end{eqnarray*}
Thus the Hilbert space projection of $f$ is 0, as it should be.
\end{example}

\begin{example} \rm
Of course $g(z) \equiv 1$ is an $H^2$ function on the disc, so its Hilbert space
projection should be $g$ itself.  Let us confirm this assertion with a calculation.

Now
\begin{eqnarray*}
Sg(z) & = & \frac{1}{2\pi} \int_{\partial D} \frac{1}{1 - z \cdot \zetabar} \, d\sigma(\zeta)  \\
                     & = & \frac{1}{2\pi i} \oint_{\partial D} \frac{\zetabar}{1 - z \cdot \zetabar} \, d\zeta \\
		     & = & \frac{1}{2\pi i} \oint_{\partial D} \frac{1}{\zeta - z} \, d\zeta = 1
\end{eqnarray*}
by the Cauchy integral formula (or simply by the calculus of residues).
\end{example}

\begin{example} \rm
Of course $h(z) \equiv z$ is an $H^2$ function on the disc, so its Hilbert space
projection should be $h$ itself.  Let us confirm this assertion with a calculation.

Now
\begin{eqnarray*}
Sh(z) & = & \frac{1}{2\pi} \int_{\partial D} \frac{\zeta}{1 - z \cdot \zetabar} \, d\sigma(\zeta) \\
                     & = & \frac{1}{2\pi i} \oint_{\partial D} \frac{1}{1 - z \cdot \zetabar} \, d\zeta \\
		     & = & \frac{1}{2\pi i} \oint_{\partial D} \frac{\zeta}{\zeta - z} \, d\zeta = z
\end{eqnarray*}
by the Cauchy integral formula (or simply by the calculus of residues).
\end{example}

\section{A Few Words About Several Complex Variables}

We record here some basic facts about the function theory of several complex variables, including a few remarks
about our three kernels in this context.  This material is presented as a setup for the application
of the Poisson-Szeg\H{o} kernel that is presented in the next section.	 We refer the reader
to [KRA1] for a more detailed treatment of several complex variables.

We work on $n$-dimensional complex Euclidean space, whose coordinates are $(z_1, \dots, z_n)$.
A function $f(z_1, \dots, z_n)$ is said to be holomorphic if it is holomorphic (in the classical
one-variable sense) in each variable separately.  This surprising definition is equivalent
to several other natural definitions of multi-variable holomorphicity:
\begin{itemize}
\item That the function have a convergent $n$-variable power series
expansion about each point in its domain;
\item That the function satisfy the Cauchy-Riemann equations in each variable
separately;
\item That the restriction of the function to each complex line $\zeta \mapsto a + b\zeta$
be holomorphic.	 That is to say, $\zeta \mapsto f(a + b\zeta)$ is holomorphic for $\zeta \in \CC$.
\end{itemize}

It is worth noting that the third of these properties implies that $\zeta
\mapsto \Re f(a + b\zeta)$ is harmonic. We call such a function {\it
pluriharmonic}. In particular, $\Re f$ is harmonic, but much more is true.
For the converse direction, if $u$ is real-valued and pluriharmonic then
$u$ is (locally) the real part of a holomorphic function.  This is proved,
as in the one-variable case, using the Poincar\'{e} lemma (see [KRA1, Ch.\ 1).
The last assertion is {\it not} true if ``pluriharmonic'' is replaced by harmonic.

In this paper we are primarily interested in the Bergman, Szeg\H{o}, and Poisson-Szeg\H{o} kernels.
These theories go through almost without change in the several variable setting.
One need only note that the Bergman space $A^2(\Omega)$, for a bounded domain $\Omega$, consists
of functions $f$ that are holomorphic on $\Omega$ and such that
$$
\int_\Omega |f(z)|^2 \, dV(z) < \infty \, ,
$$
where $dV$ is the usual Euclidean volume measure on $\Omega$.  We may stil define
$H^(\Omega)$ to be those holomorphic functions on $\Omega$ that have a harmonic majorant.  Our 
Propositions 2.1, 2.2, 2.3, 2.4, 2.5, 2.8, 2.9, 2.10 all have formalistic proofs
that are independent of dimension.  So these remain true for holomorphic functions
on $\CC^n$.  Propositions 2.6 and 2.7 require just a bit
of comment.

If $F = (f_1, \dots, f_n): \Omega \ra \Omega'$ is a holomorphic mapping of domains in $\CC^n$, 
then we may consider the $n \times n$ {\it complex Jacobian matrix}
\smallskip \\

$$
J_\CC f = \left ( \begin{array}{cccc}
              \displaystyle \frac{\partial f_1}{\partial z_1} \ \ & \ \ \displaystyle \frac{\partial f_1}{\partial z_2} \ \ & \ \  \cdots \ \ & \ \  \displaystyle \frac{\partial f_1}{\partial z_n} \\  [.19in]
              \displaystyle \frac{\partial f_2}{\partial z_1} \ \ & \ \  \displaystyle \frac{\partial f_2}{\partial z_2} \ \ & \ \  \cdots \ \ & \ \  \displaystyle \frac{\partial f_2}{\partial z_n} \\  [.19in]
						\ \ & \ \  \cdot &&   \\														       [.19in]
              \displaystyle \frac{\partial f_n}{\partial z_1} \ \ & \ \  \displaystyle \frac{\partial f_n}{\partial z_2} \ \ & \ \  \displaystyle \cdots \ \ & \ \  \displaystyle \frac{\partial f_n}{\partial z_n} \\  
		  \end{array}
	  \right ) \, .
$$
\vspace*{.135in}

\noindent This is a very natural object from the point of view of complex function theory.  But we may also write 
$f_j = u_j + i v_j$ and $z_j = x_j + i y_j$ and then consider the $2n \times 2n$ {\it real Jacobian matrix}
$$
J_\RR f = \left ( \begin{array}{ccccccc}
            \displaystyle \frac{\partial u_1}{\partial x_1} \ & \  \displaystyle \frac{\partial u_1}{\partial y_1} \ & \  \displaystyle \frac{\partial u_1}{\partial x_2} \ & \  \displaystyle \frac{\partial u_1}{\partial y_2}
	                    \ & \  \cdots \ & \  \displaystyle \frac{\partial u_1}{\partial x_n} \ & \  \displaystyle \frac{\partial u_n}{\partial y_n} \\	[.19in]
            \displaystyle  \frac{\partial v_1}{\partial x_1} \ & \  \displaystyle \frac{\partial v_1}{\partial y_1} \ & \  \displaystyle \frac{\partial v_1}{\partial x_2} \ & \  \displaystyle \frac{\partial v_1}{\partial y_2}
	                    \ & \  \cdots \ & \  \displaystyle \frac{\partial v_1}{\partial x_n} \ & \  \displaystyle \frac{\partial u_n}{\partial y_n} \\	[.19in]
			    & & & \cdots & & &				   \\ 	[.19in]
            \displaystyle  \frac{\partial u_n}{\partial x_1} \ & \  \displaystyle \frac{\partial u_n}{\partial y_1} \ & \  \displaystyle \frac{\partial u_n}{\partial x_2} \ & \  \displaystyle \frac{\partial u_n}{\partial y_2}
	                    \ & \  \cdots \ & \  \displaystyle \frac{\partial u_n}{\partial x_n} \ & \  \displaystyle \frac{\partial u_n}{\partial y_n} \\	[.19in]
            \displaystyle  \frac{\partial v_n}{\partial x_1} \ & \  \displaystyle \frac{\partial v_n}{\partial y_1} \ & \  \displaystyle \frac{\partial v_n}{\partial x_2} \ & \  \displaystyle \frac{\partial v_n}{\partial y_2}
	                    \ & \  \cdots \ & \  \displaystyle \frac{\partial v_n}{\partial x_n} \ & \  \displaystyle \frac{\partial u_n}{\partial y_n} \\
		  \end{array}
	  \right )   \, .
$$

Now it is an important fact, whose proof is just an exercise in linear algebra (see [KRA1], that these two Jacobians are related:
\smallskip \\

\noindent {\bf Proposition 2.6'} \ \ \sl Let $f: \Omega \ra \Omega'$ be a holomorphic mapping of domains
in $\CC^n$.  Then
$$
J_\RR f(z) = |J_\CC f(z)|^2 \, .
$$
\vspace*{.15in}

\rm As a consequence of this result, we can (by the same formal argument) prove this variant of Proposition 2.7:
\smallskip \\

\noindent {\bf Proposition 2.7'}  \ \ Let $\Omega_1, \Omega_2$ be domains in $\CC^n$.
Let $f: \Omega_1 \ra \Omega_2$ be biholomorphic.  Then
$$
\hbox{det}\, J_\CC f(z) \cdot K_{\Omega_2}(f(z),f(\zeta)) \cdot \overline{\hbox{det} \, J_\CC f(\zeta)} = K_{\Omega_1}(z,\zeta) \, . 
$$

The Szeg\H{o} kernel is constructed just as in the one-complex-variable context, and all the
arguments are the same.  Likewise for the Poisson-Szeg\H{o} kernel.

We conclude by noting two pieces of terminology that are special to the function theory of
several complex variables.  A domain $\Omega \ss \CC^n$ is said to be a {\it domain of holomorphy} if there
is a holomorphic function on $\O$ that cannot be analytically continued to any larger domain.\footnote{The definition
that we present here is slightly on the informal side.  For the full chapter and verse on this topic, consult [KRA1].}
A domain $\Omega \ss \CC^n$ is called {\it Levi pseudoconvex} if
\begin{enumerate}
\item[{\bf (i)}]  We can write $\Omega = \{z \in \CC^n: \rho(z) < 0\}$, where $\rho$ is $C^2$ and $\nabla \rho \ne 0$ on $\partial \Omega$;

\item[{\bf (ii)}]  For any $w \in \CC^n$ and any $P \in \partial \Omega$ such that $\sum_j [\partial \rho/\partial z_j(P)] w_j = 0$
it holds that
$$
\sum_{j,k} \frac{\partial^2 \rho}{\partial z_j \partial \zbar_k} (P) w_j \overline{w}_k \geq 0 \, .  \eqno (7.1)
$$
\end{enumerate}
The quadratic form in $(7.1)$ is called the {\it Levi form}.  A major result in the subject is that a domain
in $\CC^n$ with $C^2$ boundary is Levi pseudoconvex if and only if it is a domain of holomorphy.  We cannot
treat the matter here, but refer the reader to [KRA1] for all the details.

We conclude by taking note of this standard notation.  For $z_j = x_j + i y_j$ we have
$$
\frac{\partial}{\partial z_j} = \frac{1}{2} \left ( \frac{\partial}{\partial x_j} - i \frac{\partial}{\partial y_j} \right ) \qquad \hbox{and} \qquad
  \frac{\partial}{\partial \zbar_j} = \frac{1}{2} \left ( \frac{\partial}{\partial x_j} + i \frac{\partial}{\partial y_j} \right ) \, .
$$
We note that
\begin{eqnarray*}
\frac{\partial}{\partial z_j} \, z_j = 1 \quad & \quad & \quad \frac{\partial}{\partial \zbar_j} \, z_j = 0 \\ 
\frac{\partial}{\partial \zbar_j} \, z_j = 0 \quad & \quad & \quad \frac{\partial}{\partial \zbar_j} \, \zbar_j = 1  \, .
\end{eqnarray*}

Likewise, if we set $dz_j = dx_j + i dy_j$ and $d\zbar_j = dx_j - i dy_j$ then 
$\langle dz_j, \partial/\partial z_j \rangle = 1$, $\langle d\zbar_j, \partial/\partial \zbar_j \rangle = 1$,
and all other pairings are equal to 0.

\section{More on the Poisson-Szeg\H{o} Kernel}

If $g = (g_{jk})$ is a Riemannian metric on a domain $\Omega$
in Euclidean space, then there is a second-order partial
differential operator, known as the {\it Laplace-Beltrami
operator}, that is invariant under isometries of the metric.
In fact, if $g$ denotes the determinant of the metric matrix $g$, and if
$(g^{jk})$ denotes the inverse matrix, then
this partial differential operator is defined to be
$$
{\cal L} = \frac{2}{g} \sum_{j,k} \left \{ \frac{\partial}{\partial \bar{z}_j} \left ( g g^{jk} \frac{\partial}{\partial z_k} \right )
                + \frac{\partial}{\partial z_k} \left ( g g^{jk} \frac{\partial}{\partial \bar{z}_k} \right ) \right \} \, .
$$

Now of course we are interested in artifacts of the Bergman theory.  If $\Omega \ss \CC^n$ is a bounded
domain and $K= K_\Omega$ its Bergman kernel, then Proposition 2.8 tells us that
$K(z,z) > 0$ for all $z \in \Omega$.  Then it makes sense to define
$$
g_{jk}(z) = \frac{\partial^2}{\partial z_k \partial \zbar_k} \log K(z,z)
$$
for $j, k = 1, \dots, n$.
Then Proposition 2.7 can be used to demonstrate that this metric---which is in fact a K\"{a}hler metric on $\Omega$---is 
invariant under biholomorphic mappings of $\Omega$.  In other words, any biholomorphic
$\Phi: \Omega \ra \Omega$ is an isometry in the metric $g$.  This is the celebrated {\it Bergman metric}.

If $\Omega \ss \CC^n$ is the unit ball $B$, then the Bergman kernel is given
by
$$
K_B(z, \zeta) = \frac{1}{V(B)} \cdot \frac{1}{(1 - z\cdot \zetabar)^{n+1}} \, ,
$$
where $V(B)$ denotes the Euclidean volume of the domain $B$ (this by a calculation similar to,
more complicated than, that in Example 3.1---see [KRA1] for the details).
Then
$$
\log K(z,z) = - \log V(B) - (n+1) \log (1 - |z|^2) .
$$
Further,
$$
\frac{\partial}{\partial z_j} \bigl ( - (n+1) \log (1 - |z|^2) \bigr )  = 
                                         (n+1) \frac{\bar{z}_j}{1 - |z|^2} 
$$
and
\begin{eqnarray*}
 \frac{\partial^2}{\partial z_j \partial \bar{z}_k} \bigl ( - (n+1) \log (1 - |z|^2 \bigr ) &  = &
                                                  (n+1) \left [ \frac{\delta_{jk}}{1 - |z|^2} + \frac{\bar{z}_j z_k}{(1 - |z|^2)^2} \right ] \\
              & = & \frac{(n+1)}{(1 - |z|^2)^2} \bigl [ \delta_{jk} (1 - |z|^2) + \bar{z}_j {z}_k \bigr ] \\
           & \equiv & g_{jk}(z) .
\end{eqnarray*}

When $n=2$  we have
$$
g_{jk}(z) = \frac{3}{(1 - |z|^2)^2} \bigl [ \delta_{jk}(1 - |z|^2) + \bar{z}_j z_k \bigr ] .
$$
Thus
$$
  \bigl ( g_{jk}(z) \bigr ) = \frac{3}{(1 - |z|^2)^2}
             \left ( \begin{array}{lr}
                 1 - |z_2|^2 & \bar{z}_1 z_2 \\
                 \bar{z}_2 z_1 & 1 - |z_1|^2
                     \end{array}
             \right ) .
$$
Let 
$$
\biggl ( g^{jk}(z) \biggr )_{j,k=1}^n
$$
represents the inverse of the matrix 
$$
\biggl ( g_{jk}(z) \biggr )_{j,k=1}^n \ \ .
$$
Then an elementary computation shows that
$$
\biggl ( g^{jk}(z) \biggr )_{j,k=1}^n
     = \frac{1 - |z|^2}{3} \left ( \begin{array}{lr}
                         1  - |z_1|^2 & - z_2 \bar{z}_1 \\
                        - z_1 \bar{z}_2 & 1 - |z_2|^2 
                                   \end{array}
                           \right )
     = \frac{1 - |z|^2}{3} \bigl ( \delta_{jk} - \bar{z}_j z_k \bigr )_{j,k} .
$$
Let
$$
g \equiv \det \biggl ( g_{jk}(z) \biggr ) .
$$
Then
$$
 g = \frac{9}{(1 - |z|^2)^3} . 
$$

Now let us calculate.
If $\bigl ( g_{jk}\bigr )^n_{j,k=1}$ is the Bergman metric on the
ball in $\CC^n$ then we have
$$
\sum_{j,k} \frac{\partial}{\partial \bar{z}_j} \bigl ( g g^{jk} \bigr ) = 0
$$
and 
$$ 
\sum_{j,k} \frac{\partial}{\partial z_j} \bigl ( g g^{jk} \bigr ) = 0 .
$$
We verify these assertions in detail in dimension $2:$
Now
\begin{eqnarray*}
 g g^{jk} & = & \frac{9}{(1 - |z|^2)^3} \cdot \frac{1 - |z|^2}{3} (\delta_{jk} - \bar{z}_j z_k ) \\
          & = & \frac{3}{(1- |z|^2)^2} (\delta_{jk} - \bar{z}_j z_k ) .
\end{eqnarray*}
It follows that
$$
\frac{\partial}{\partial \bar{z}_j} \biggl [ g g^{jk} \biggr ] =
              \frac{6 z_j}{(1 - |z|^2)^3} \bigl ( \delta_{jk} - \bar{z}_j z_k \bigr )      
                                    - \frac{3z_k}{(1 - |z|^2)^2} .
$$
Therefore
\begin{eqnarray*}
          \sum_{j,k = 1}^2 \frac{\partial}{\partial \bar{z}_j} \biggl [ g g^{jk} \biggr ] & = & 
                   \sum_{j,k = 1}^2 \left [ \frac{6 z_j (\delta_{jk} - \bar{z}_j z_k)}{(1 - |z|^2)^3}
                             - \frac{3 z_j}{(1 - |z|^2)^2} \right ] \\
               & = & 6 \sum_k \frac{z_k}{(1 - |z|^2)^3} -
                     6 \sum_{j,k} \frac{|z_j|^2 z_k}{(1 - |z|^2)^3} 
                     - 6 \sum_k \frac{z_k}{(1 - |z|^2)^2}  \\
               & = &  6 \sum_{j} \frac{z_k}{(1 - |z|^2)^2}
                       - 6 \sum_k \frac{z_k}{(1 - |z|^2)^2} \\
               & = & 0 .
\end{eqnarray*}
The other derivative is calculated similarly.

Our calculations show that, on the ball in $\CC^2,$
\begin{eqnarray*}
{\cal L} & \equiv & \frac{2}{g} \sum_{j,k} \left \{ \frac{\partial}{\partial \bar{z}_j}
                                      \left ( g g^{jk} \frac{\partial}{\partial z_k} \right ) +
                                      \frac{\partial}{\partial z_k}
                                      \left ( g g^{jk} \frac{\partial}{\partial \bar{z}_j} \right ) \right \} \\
                 & = & 4 \sum_{j,k} g^{jk} \frac{\partial}{\partial \bar{z}_j} \frac{\partial}{\partial z_k} \\
                 & = & 4 \sum_{j,k} \frac{1 - |z|^2}{3} \bigl ( \delta_{jk} - \bar{z}_j z_k \bigr )
                           \frac{\partial^2}{\partial z_k \partial \bar{z}_j} .
\end{eqnarray*}

Now the interesting fact for us is encapsulated in the following proposition:

\begin{proposition} \sl
The Poisson-Szeg\H{o} kernel on the ball $B$ solves the Dirichlet problem
for the invariant Laplacian ${\cal L}$.  That is to say, if $f$ is a continuous
function on $\partial B$ then the function
$$
u(z) = \left \{ \begin{array}{lcr}
             \int_{\partial B} {\cal P}(z, \zeta) \cdot f(\zeta) \, d\sigma(\zeta) & \hbox{if} & z \in B \\  [.15in]
	       f(z)                                                                & \hbox{if} & z \in \partial B 
		\end{array}
       \right.
$$
is continuous on $\overline{B}$ and is annihilated by ${\cal L}$ on $B$.
\end{proposition}

This fact is of more than passing interest.  In one complex variable, the study of holomorphic
functions on the disc and the study of harmonic functions on the disc are inextricably linked
because the real part of a holomorphic function is harmonic and conversely.  Such is not
the case in several complex variables.  Certainly the real part of a holomorphic function
is harmonic.  But in fact it is more:  such a function is {\it pluriharmonic}.  For the
converse direction, any real-valued pluriharmonic function is locally the real part of a holomorphic
function.  This assertion is false if ``pluriharmonic'' is replaced by ``harmonic''.  

Thus, in some respects, it is inappropriate to study holomorphic functions on the ball in $\CC^n$ using
the Poisson kernel.  In view of Proposition 8.1, the Poisson-Szeg\H{o} kernel is much more
apposite.  As an instance, Adam Koranyi [KOR] made decisive use of this observation
in his study of the boundary behavior of $H^2(B)$ functions.

We shall spend the rest of this section discussing the proof of the proposition.  There are two things
to show:
\begin{enumerate}
\item[{\bf (1)}]  That the function $u$ is annihilated by the invariant Laplacian ${\cal L}$ on $B$;
\item[{\bf (2)}]  That the function $u$, defined piecewise in the proposition, is actually continuous
on $\overline{B}$.
\end{enumerate}
Most of our efforts will be used to dispatch {\bf (1)}.  We shall make a few remarks at the end about {\bf (2)}.

In fact there are a number of known methods for addressing {\bf (1)}.  E. M. Stein, in [STE], verified that
in fact the function $u$ satisfies a suitable mean value condition (as in harmonic function theory) on suitable
balls in $B$.  Invoking a theorem of Godement, he was then able to conclude that $u$ was ${\cal L}$-harmonic.
L. Hua, in [HUA], gave an alternative argument.  Yet another approach may be found in [HEL].  We take
here a more elementary approach and actually {\it calculate} ${\cal L}u$.

Now
$$
{\cal L} u = {\cal L} \int_{\partial B} {\cal P}(z, \zeta) \cdot f(\zeta) \, d\sigma(\zeta)
            = \int_{\partial B} \biggl [ {\cal L}_z {\cal P} (z,\zeta)\biggr  ] \cdot f(\zeta) \, d\sigma(\zeta) \, .
$$
Thus it behooves us to calculate ${\cal L}_z {\cal P}(z, \zeta)$.   Now we shall calculate this
quantity for each fixed $\zeta$.  Thus, without loss of generality, we may compose with a unitary
rotation and suppose that $\zeta = (1 + i0, 0 + i0)$ so that (in complex dimension 2)
$$								   
\P = c_2 \cdot \frac{(1 - |z|^2)^2}{|1 - z_1|^4} \, .
$$
This will make our calculations considerably easier.
						   
By brute force, we find that

\small

\begin{align}
\frac{\partial \P}{\partial \zbar_1} & = -2(1 - z_1)(1 - |z|^2) \cdot \left [ \frac{-1 + z_1 + |z_2|^2}{|1 - z_1|^6} \right ] \notag \\ 
\frac{\partial^2 \P}{\partial \zbar_1 \partial z_1} & = \frac{-2}{|1 - z_1|^6} \cdot \left [ - |z_1|^2 - |z_1|^2|z_2|^2 + 3|z_2|^2 - z_1 |z_2|^2 \right. \notag \\ 
             & \null \ \ \ \ \quad \left. - 2|z_2|^4 - 1 + z_1 + \zbar_1 - \zbar_1|z_2|^2 \right ] \notag \\  
\frac{\partial^2 \P}{\partial \zbar_1 \partial z_2} & = \frac{-2(1 - z_1)}{|1 - z_1|^6} \cdot \left [ 2\zbar_2 
        - \zbar_2 z_1 - 2 \zbar_2 |z_2|^2 - \zbar_2 |z_1|^2 \right ] \notag \\  
\frac{\partial^2 \P}{\partial z_1 \partial \zbar_2} & = \frac{-2(1 - \zbar_1)}{|1 - z_1|^6} \cdot \left [ 2z_2 
        - z_2 \zbar_1 - 2 z_2 |z_2|^2 - z_2 |z_1|^2 \right ] \notag \\	 
\frac{\partial \P}{\partial z_2} & = \frac{-2z_2 + 2|z_1|^2 z_2 + 2|z_2|^2 z_2}{|1 - z_1|^4} \notag \\ 
\frac{\partial^2 \P}{\partial z_2 \partial \zbar_2} & = \frac{-2 + 2|z_1|^2 + 4|z_2|^2}{|1 - z_1|^4} \tag{8.1} \\  \notag
\end{align}

\normalsize

Now we know that, in complex dimension two,

\small

\begin{eqnarray*}
{\cal L}_z \P(z, \zeta) & = & \frac{4}{3}(1 - |z|^2) \cdot (1 - |z_1|^2) \cdot \frac{\partial^2 \P_z}{\partial z_1 \partial \zbar_1} +
	   \frac{4}{3}(1 - |z|^2) \cdot (- \zbar_1 z_2) \cdot \frac{\partial^2 \P_z}{\partial z_2 \partial \zbar_1} \\
	   && \quad + \frac{4}{3}(1 - |z|^2) \cdot (- \zbar_2 z_1) \cdot \frac{\partial^2 \P_z}{\partial z_1 \partial \zbar_2} +
	   \frac{4}{3}(1 - |z|^2) \cdot (1 - |z_2|^2) \cdot \frac{\partial^2 \P_z}{\partial z_2 \partial \zbar_2} \, .
\end{eqnarray*}

\normalsize

\noindent Plugging the values from (8.1) into this last equation gives

\small

\begin{eqnarray*}
{\cal L}_z \P(z, \zeta) & = & \frac{4}{3}(1 - |z|^2) \cdot (1 - |z_1|^2) \cdot 
	 \frac{-2}{|1 - z_1|^6} \cdot \biggl [ - |z_1|^2 - |z_1|^2|z_2|^2 \\
	  && \null \ \ \ \ \quad + 3|z_2|^2 - z_1 |z_2|^2 
              - 2|z_2|^4 - 1 + z_1 + \zbar_1 - \zbar_1|z_2|^2 \biggr ]  \\  
	  && \ + \frac{4}{3}(1 - |z|^2) \cdot (- \zbar_1 z_2) \\ 
	  && \quad \ \ \ \times  \frac{-2(1 - z_1)}{|1 - z_1|^6} \cdot \biggl [ 2\zbar_2 
        - \zbar_2 z_1 - 2 \zbar_2 |z_2|^2 - \zbar_2 |z_1|^2 \biggr ]  \\  
	  && \ + \frac{4}{3}(1 - |z|^2) \cdot (- \zbar_2 z_1) \\ 
	  && \quad \ \ \ \times  \frac{-2(1 - \zbar_1)}{|1 - z_1|^6} \cdot \biggl [ 2z_2 
        - z_2 \zbar_1 - 2 z_2 |z_2|^2 - z_2 |z_1|^2 \biggr ]  \\
	  && \ + \frac{4}{3}(1 - |z|^2) \cdot (1 - |z_2|^2) \cdot |1 - z_1|^2 \cdot
	    \frac{-2 + 2|z_1|^2 + 4|z_2|^2}{|1 - z_1|^6} \, . \\
\end{eqnarray*}

\def\z1{z_1}
\def\z1bar{\overline{z}_1}
\def\z2{z_2}
\def\z2bar{\overline{z}_2}

\normalsize

Multiplying out the terms, we find that

\small

\begin{eqnarray*}
{\cal L}_z \P(z, \zeta) & = & \frac{-2}{|1 - z_1|^6} \cdot \biggl [ - |z_1|^2  - 4|z_1|^2|z_2|^2 + 3|z_2|^2 - z_1 |z_2|^2 - 2|z_2|^4 - 1 \\
			&& \quad + z_1 + \zbar_1 - \zbar_1|z_2|^2 + |z_1|^4 + |z_1|^4|z_2|^2 + z_1 |z_1|^2|z_2|^2 \\
			&& \quad + 2|z_1|^2 |z_2|^4 + |z_1|^2 - z_1|z_1|^2 - \zbar_1|z_1|^2 + \zbar_1|z_1|^2|z_2|^2 \biggr ] \\
			&& \quad - \frac{2}{|1 - z_1|^6} \cdot \biggl [ -2 \zbar_1 |z_2|^2 + 3|z_1|^2 |z_2|^2 + 2|z_2|^4 \zbar_1 + \zbar_1 |z_2|^2 |z_1|^2 \\
			&& \quad - z_1 |z_1|^2 |z_2|^2 - 2 |z_1|^2 |z_2|^4 - |z_2|^2 |z_1|^4 \biggr ] \\
			&& \quad - \frac{2}{|1 - z_1|^6} \cdot \biggl [ -2 z_1 |z_2|^2 + 3|z_1|^2 |z_2|^2 + 2|z_2|^4 z_1 + z_1 |z_2|^2 |z_1|^2 \\
			&& \quad - \zbar_1 |z_1|^2 |z_2|^2 - 2 |z_1|^2 |z_2|^4 - |z_2|^2 |z_1|^4 \biggr ] \\
			&& \quad - \frac{2}{|1 - z_1|^6} \cdot \biggl [ 1 - |z_1|^2 - 3|z_2|^2 + |z_1|^2 |z_2|^2 + 2|z_2|^4 - z_1 + z_1|z_1|^2 \\
			&& \qquad + 3 z_1 |z_2|^2 - z_1 |z_1|^2 |z_2|^2 - 2 z_1 |z_2|^4 - \zbar_1 + \zbar_1 |z_1|^2 + 3 \zbar_1 |z_2|^2 - \zbar_1 |z_1|^2 |z_2|^2 \\
			&& \qquad - 2 \zbar_1 |z_2|^4 + |z_1|^2 - |z_1|^4 - 3|z_1|^2 |z_2|^2 + |z_1|^4 |z_2|^2 + 2|z_1|^2 |z_2|^4 \biggr ] \, .	   \\
\end{eqnarray*}

\normalsize

And now if we combine all the terms in brackets a small miracle happens:  everything cancels.  The result is
$$
\L_z \P(z, \zeta) \equiv 0.
$$

That is half of our task.  For the other half, notice that the Poisson-Szeg\H{o} kernel
$$
\P(z, \zeta) = c_n \cdot \frac{(1 - |z|^2)^2}{|1 - z\cdot \zeta|^4}
$$
has these properties:
\begin{enumerate}
\item[{\bf (2)}]  $\P(z,\zeta) > 0$ for all $z \in B$, $\zeta \in \partial B$;

\item[{\bf (1)}]  $\int_{\partial B} \P(z,\zeta) \, d\sigma(\zeta) = 1 \qquad \forall z \in B$;

\item[{\bf (3)}]  If $\epsilon > 0$, then $\lim_{|z| \ra 1} \P(z, \zeta) = 0,$
          uniformly for $|z - \zeta| \geq \epsilon.$  
\end{enumerate}
These properties are clear by inspection.  They are analogous to the properties
of a ``standard family of kernels'' as discussed in [KRA5] and [KAT].  And it is a general
and well-known fact that, with these desiderata, the integral
$$
\int_{\partial B} \P(z, \zeta) \, f(\zeta) \, d\sigma(\zeta)
$$
will converge back to $f(P)$ as $B \ni z \ra P \in \partial B$.  That is part
{\bf (2)} of what we wished to prove.

\section{Concluding Remarks}

The Bergman, Szeg\H{o}, and Poisson-Szeg\H{o} kernels are today
central to the study of harmonic analysis in the complex
domain. Although we did not treat the matter here, the
Szeg\H{o} kernel is usually a classical singular integral
kernel on the boundary $\partial \Omega$, hence the
Cald\'{e}ron-Zygmund theory may be brought to bear on the
subject (see [KRA6] for more on this particular topic)---see
particularly its more general formulation in [COW].  The
Bergman kernel is, by contrast, a classical Hilbert integral
and may be studied using slightly different tools (see [PHS]).

The three integral operators featured here are crucial to the regularity theory
of partial differential equations---particularly the inhomogenous Cauchy-Riemann equations
and the Laplacian---see [KRA4].   They help us to understand function theory, the
construction of holomorphic function with particular properties, and the nature of
domains of holomorphy.  They are important tools.

In the present paper we have only given a light introduction to, and an invitation
to, the study of the kernels of Bergman, Poisson, and Szeg\H{o}.  We can only
hope that the reader is tempted to explore further.
\newpage

\null \vspace*{.5in}

\section*{\sc References}
\vspace*{.2in}

\begin{enumerate}

\item[{\bf [APF]}] L. Apfel, Localization properties and
boundary behavior of the Bergman kernel, thesis, Washington
University in St. Louis, 2003.

\item[{\bf [ARO]}]  N. Aronszajn, Theory of reproducing kernels, {\it Trans.\
Am.\ Math.\ Soc.} 68(1950), 337-404. 

\item[{\bf [BER]}]  S. Bergman, {\it The Kernel Function and Conformal Mapping},
Am. Math. Soc., Providence, RI, 1970.
									      
\item[{\bf [BLK]}]  B. Blank and S. G. Krantz, {\it Calculus:  Multivariable},
Key College Press, Emeryville, CA, 2006.  

\item[{\bf [COW]}] R. R. Coifman and G. L. Weiss, {\it Analyse
Harmonique Non-Commutative sur Certains Espaces Homogenes},
Springer Lecture Notes vol. 242, Springer Verlag, Berlin,
1971.

\item[{\bf [GK]}] R. E. Greene and S. G. Krantz, {\it Function
Theory of One Complex Variable}, $2^{\rm nd}$ ed., American
Mathematical Society, Providence, RI, 2002.

\item[{\bf [HEL]}]  S. Helgason, {\it Differential Geometry and Symmetric Spaces},
Academic Press, New York, 1962.

\item[{\bf [HUA]}]  L.-K. Hua, {\it Harmonic Analysis of
Functions of Several Complex Variables in the Classical
Domains}, American Mathematical Society, Providence, 1963.

\item[{\bf [KAT]}] Y. Katznelson {\it An Introduction to
Harmonic Analysis}, Wiley, New York, 1968.
				       
\item[{\bf [KOR]}]  A. Koranyi, The Poisson integral for generalized half-planes
and bounded symmetric domains, {\it Annals of Math.} 82(1965), 332--350.
	
\item[{\bf [KRA1]}]  S. G. Krantz, {\it Function Theory of
Several Complex Variables}, $2^{\rm nd}$ ed., American
Mathematical Society, Providence, RI, 2001.

\item[{\bf [KRA2]}]  S. G. Krantz, {\it Cornerstones of
Geometric Function Theory: Explorations in Complex Analysis},
Birkh\"{a}user Publishing, Boston, 2006.

\item[{\bf [KRA3]}]  S. G. Krantz, Calculation and estimation of the Poisson kernel,
{\it J. Math.\ Anal.\ Appl.} 302(2005)143--148. 

\item[{\bf [KRA4]}]  S. G. Krantz, {\it Partial Differential Equations and Complex
Analysis}, CRC Press, Boca Raton, FL, 1992.

\item[{\bf [KRA5]}]  S. G. Krantz, {\it A Panorama of Harmonic Analysis},
Mathematical Association of America, Washington, D.C., 1999.

\item[{\bf [KRA6]}]  S. G. Krantz, {\it Explorations in Harmonic Analysis 
with Applications to Complex Function Theory and the Heisenberg Group},
Birkh\"{a}user Publishing, Boston, MA, to appear.

\item[{\bf [NRSW]}] A. Nagel, J.-P. Rosay, E. M. Stein, and S.
Wainger, Estimates for the Bergman and Szeg\"{o} kernels in
$\CC^2,$ {\em Annals Math.} 129(1989), 113-149.

\item[{\bf [PHS]}] D. H. Phong and E. M. Stein, Hilbert integrals, singular
integrals, and Radon transforms. I. {\it Acta Math.} 157(1986), 99--157.

\item[{\bf [STE]}] E. M. Stein, {\it Boundary Behavior of
Holomorphic Functions of Several Complex Variables}, Princeton
University Press, Princeton, 1972.

\end{enumerate}

\end{document}